\documentclass[12pt]{article}
\usepackage{times}
\usepackage{amsmath}
\usepackage{amssymb}
\usepackage{color}
\newcommand{\secref}[1]{Section~\ref{#1}}

\newcommand{\bsp}{B} 
\newcommand{\sft}{\gamma}
\newcommand{\dg}{\delta}
\newcommand{\be}{\mathbf{e}}
\newcommand{\pdg}{p} 
\newcommand{\convol}{\star}
\newcommand{\trns}{\mathtt{t}}
\newtheorem{corollary}{Corollary}[section]
\newtheorem{theorem}{Theorem}[section]
\newcommand{\lst}{n}
\newcommand{\itoj}{i:j}
\newcommand{\gs}{\sigma}
\newcommand{\knts}{\sft:\sft+d+1}  
\newcommand{\tsft}{t_{\knts}}  
\newcommand{\divdif}{\mathord
   {\kern.43em\vrule width.6pt height5.6pt depth-.28pt\kern-.43em\Delta}}

\newcommand{\mc}{\mathbf{c}}
\newcommand{\R}{\mathbb{R}}
\newenvironment{proof}{
    {\bf Proof}\quad}{\hbox{\ }\ \hfil $|||$\break 
            }

\title{Computing SIAC spline coefficients} 
\author{J\"org Peters} 

\begin{document}
\maketitle
\begin{abstract}
The Discontinuous Galerkin (DG) method
applied to hyperbolic differential equations outputs
weakly-linked polynomial pieces.
Post-processing these pieces by 
Smoothness-Increasing Accuracy-Conserving (SIAC) convolution 
with B-splines can improve the accuracy of the output and
yield superconvergence.
SIAC convolution is considered optimal if the 
SIAC kernels, in the form of a linear combinations of B-splines
of degree $d$, reproduce polynomials of degree $2d$.
This paper derives simple formulas for computing the optimal
SIAC spline coefficients.
\end{abstract}
\section{Introduction}
The Discontinuous Galerkin (DG) method is widely used to solve
hyperbolic partial differential equations (PDEs). The lack of continuity 
between the elements is appropriate for modeling the
weak flux constraints between elements and computationally convenient:
the discontinuity allows for a flexible discretization of the PDE
that locally adjusts the polynomial degree and element spacing,
and the discontinuity increases opportunities for parallelism when
stepping forward in a simulation.
However, except near jump discontinuities,  the inter-element discontinuities 
do not agree with the smoothness of the expected outcome. 

Filtering, in  particular Smoothness-Increasing Accuracy-Conserving (SIAC)
filtering, has
been proposed to smoothly connect elements while
maintaining the order of the accuracy of the original DG solution.
Remarkably, such post-filtering by convolution applied to approximate DG
solutions is not only convenient for downstream applications
such as stream line tracing \cite{journals/jscic/WalfischRKH09}, 
but improves the accuracy of the 
resulting output as a solution to hyperbolic partial differential equations.
Already Bramble and Schatz \cite{Bramble:1977:HOL} showed that,
for a wide class of elliptic boundary value problems and uniform subdivision
of the domain, averaging of the output can yield superconvergence,
i.e.\ a more accurate approximation to the solution than the degree of the 
elements suggests.
Underlying superconvergence is the fact that certain integral norms, 
called moment norms \cite{MockLax} or negative-order norms,
converge faster than expected and can bound the error 
of the convolved output.
In the context of linear hyperbolic PDEs this fact was 
convincingly demonstrated in \cite{journals/moc/CockburnLSS03}.


Starting with \cite{journals/siamsc/CurtisKRS07},
a series of papers has generalized SIAC filtering of DG output
from the prototypical case of
linear equations with periodic boundary conditions over a uniform mesh
to non-uniform scenarios and 
spatial dimensions two and three, including
structured and unstructured bivariate and tetrahedral meshes
\cite{journals/siamnum/MirzaeeJRK11,journals/siamsc/MirzaeeKRK13,%
journals/jscic/MirzaeeRK14}.

SIAC filtering was
authoritatively reviewed during the 2014 Icosahom conference,
notably in the planary talk by Jennifer Ryan.
Newest results were presented in
the minisymposium on post-processing DG
solutions \cite{mirz:2014}, organized by Mirzagar and Ryan.
Among the recent advances were simple formulas 
that allow solving for the optimal coefficients of uniform
SIAC B-spline  filters (note by contrast that \cite{journals/jscic/MirzaeeRK12} 
used Gaussian quadrature too determine the entries of the 
corresponding constraint matrix).
Since the SIAC approach has been extended to non-uniform meshes, 
formulas corresponding to the general case of B-splines with non-uniform knot
spacing can improve computational efficiency.
The following pages develop such formulas.

\secref{sec:basics} succinctly reviews B-splines, SIAC filtering and convolution
to the extent needed for the result.
\secref{sec:formula} derives the entries of the 
constraint matrix whose solution yields the 
optimal coefficients for filters that post-process DG solutions
with splines over non-uniform knot sequences.



\section{Convolution and B-splines}
\label{sec:basics}
A succint but comprehensive treatment of 
B-splines can be found in Carl de Boor's summary \cite{Boor:2002:BB}
(see also \cite{Schumaker:1981:SF}).
There are a number of ways to derive or define B-splines,
for example as the smoothest class of piecewise polynomials over
a given support. The subclass of uniform B-splines additionally admits
definitions via convolution that are efficiently carried out in Fourier space.
That definition is handy when deriving optimal coefficients of uniform
SIAC B-spline filters.
For our purpose, a more convenient definition is via 
divided differences $\divdif(t_{i:j})h$%
\footnote{The authoritative
survey \cite{deboor} advertises the symbol $\divdif$ for divided differences 
over alternatives such as $[t_{i:j}]h$ or $h[t_{i:j}]$.}.
The notation $i:j$ is short for the sequence $i,i+1,\ldots,j-1,j$
and $t_{i:j}$ stands correspondingly for the sequence of real numbers
$t_i,t_{i+1},\ldots,t_j$.
For a sufficiently smooth univariate real-valued
function $h$ with $k$th derivative $h^{(k)}$,
divided differences are defined by
\begin{align}
   \divdif(t_i) h &:= h(t_i),
   \hskip0.25\linewidth
   \text{ and for }
   j > i 
   \notag
   \\
   \divdif(t_{i:j}) h &:= 
   \begin{cases} 
      (\divdif(t_{i+1:j}) h - \divdif(t_{i:j-1}) h)/ (t_j-t_i), & 
      \text{ if } t_i\ne t_j,\\
      \frac{h^{(j-i)}}{(j-i)!} (t_i), &
      \text{ if } t_i = t_j.
   \end{cases}
    \label{eq:divdif}
\end{align}
If $t_{\itoj}$ is a non-decreasing sequence, we call its elements
$t_\ell$ \emph{knots} 
and define the \emph{B-spline with knot sequence $t_{\itoj}$} as
\begin{align}
   \bsp(x|t_{\itoj}) &:= (t_{j}-t_i)\ \divdif(t_{\itoj})
   (\max\{(\cdot-x),0\})^{d},
   \label{eq:defa}
\end{align}
where $\divdif(t_{\itoj})$ acts on the function 
$h: t \to (\max\{(t-x),0\})^{d}$ for a given $x \in \R$.
A B-spline is a non-negative piecewise polynomial of degree $d$
with support on the interval $[t_i..t_j)$.
If $\mu$ is the multiplicity of the 
number $t_\ell$ in the sequence $t_{\itoj}$,
then $\bsp(x|t_{\itoj})$ is at least $d-\mu$ times continuously differentiable
at $t_\ell$.

The goal of SIAC filtering is to smooth out the transitions between
polynomial pieces on consecutive intervals
$\pdg_j: [t_j..t_{j+1})\to \R$, $j=0..\lst$ that are output by DG computations
and 
typically do not join continuously.
To this end, we convolve the piecewise output
with a linear combination of B-splines. 
The convolution $f\convol g$
of  
a function $f$ with a function $g$ is
defined as
\begin{align}
   (f\convol g)(x) := \int_\R f(t) g (x-t) dt,
   \label{eq:convoldef}
\end{align}
for every $x$ where the integral exists.
When $g \ge 0$ and $\int_\R g =1$ then
the convolution has special,  desirable properties:
if $f$ is non-negative, (directionally) monotone or convex
then so is $f\convol g$.
Moreover, the graph of $f\convol g$ is in the convex hull 
of the graph of $f$.
Convolution is commutative, associative and distributive.

For convolving splines with functions $g := \pdg_j$,
we make use of Peano's formula: 
\begin{align}
   \frac{1}{k!} \int_\R
   \bsp(t|t_{0:k}) 
   g^{(k)}(t) dt
   =
   \divdif(t_{0:k}) g.  
   \label{eq:peano}
\end{align}
Here $g^{(k)}$ denotes the $k$th derivative of the univariate function $g$.
To be able to interpret Peano's formula as a convolution formula for monomials,
we select $g$ so that $g^{(k)}(t) = (k!)(t-x)^\dg$.
The choice $g(t) := \binom{k+\dg}{\dg}^{-1} (t-x)^{k+\dg}$
accomplishes this. Then \eqref{eq:peano} implies
for the alternating monomial $(-\cdot)^\dg: t \to (-t)^\dg$,
\begin{align}
   \left(\bsp(\cdot|t_{0:k}) \convol (-\cdot)^\dg\right) (x)
   &=^\eqref{eq:convoldef}
   \int_\R \bsp(t|t_{0:k}) (-(x-t))^\dg dt
   =
   \frac{1}{k!} \int_\R \bsp(t|t_{0:k}) g^{(k)}(t) dt
   \label{eq:peanopol}
   \\
   \notag
   &=^\eqref{eq:peano}
   \binom{k+\dg}{k}^{-1}
   \divdif{t_{0:k}} (t-x)^{k+\dg}.
\end{align}
(The divided difference in the last expression applies 
to the variable $t$. Therefore $\divdif{t_{0:k}} (t-x)^{k+\dg}$
does not depend on $t$, but only on the sequence $t_{0:k}$.)

\section{Optimal convolution coefficients}
\label{sec:formula}
A spline SIAC convolution kernel $K: \R\to\R$ is a piecewise polynomial of
degree $d$. The function $K$ is considered optimal if 
\begin{align}
   (K \convol (\cdot)^\dg) (x) = x^\dg, \qquad \dg = 0,\ldots,2d,
   \label{eq:convol}
\end{align}
i.e.\ if  convolution of $K$ with monomials reproduces the monomials 
up to the maximal degree $2d$.
The choice of interest for SIAC convolution on the interval $[t_0..t_{d+1})$ 
is $K(x) := \sum^d_{\sft=-d}  c_\sft B(x|\tsft)$,
a spline of degree $d$ with leftmost knot $t_{-d}$.
To satisfy the polynomial equations \eqref{eq:convol},
we want to determine the coefficients $c_{-d},\ldots,c_d$ so that 
\begin{align}
   &\Bigl( \sum^d_{\sft=-d}  c_\sft B(\cdot|\tsft)
   \convol (-\cdot)^\dg \Bigr) (x)
   =
   (-x)^\dg, \qquad \dg = 0,\ldots,2d.
   \tag{\ref{eq:convol}'}
\end{align}
\begin{theorem}
The vector of optimal SIAC convolution coefficients 
$\mc := [c_{-d},\ldots,c_d]^\trns \in \R^{2d+1}$ is
\begin{align}
   \mc &:= M_0^{-1} \be_1,
   \qquad
   \be_i(\dg) := 
   \begin{cases}
      1 & \text{ if } \dg=i\\
      0 & \text{ else },
   \end{cases}
   \label{eq:M}
   \\
   \notag
   &M_0 := 
   \left[ \begin{matrix}
   \divdif{\tsft} t^{d+1+\dg}
   \end{matrix} \right]_{\dg=0:2d, \sft=-d:d}.
\end{align}
\end{theorem}
The matrix $M_0$ is of size $(2d+1)\times(2d+1)$. Its
entries $M_0(k,\ell)\in \R$ in row $k$ and column $\ell$
are defined by $\dg = k-1$ and $\sft=\ell-d-1$ and 
do not depend on $t$, but on $\tsft$.

\begin{proof}
By \eqref{eq:peanopol}, the system of polynomial equations (\ref{eq:convol}')
in $x$ is equivalent to 
\begin{align}
   M(x)\mc &:= 
   \left[ \begin{matrix}
   \divdif{\tsft} (t-x)^{d+1+\dg}]
   \end{matrix} \right]_{\dg=0:2d, \sft=-d:d}
   \mc 
   \\
   \notag
   &=
   \binom{d+1+\dg}{\dg}
   \left[ \begin{matrix}
   (-x)^\dg
   \end{matrix} \right]_{\dg=0:2d}.
   \tag{\ref{eq:convol}''}
\end{align}
Since (\ref{eq:convol}'') has to hold for $x=0$,
setting $x=0$ in (\ref{eq:convol}'') 
yields the following $2d+1\times 2d+1$ system of equations
\begin{equation}
   \left[ \begin{matrix}
   \binom{d+1+\dg}{\dg}^{-1}
   \divdif{\tsft} t^{d+1+\dg} 
   \end{matrix} \right]_{\dg=0:2d, \sft = -d:d} 
   \mc
   =
   \be_1.
   \tag{\ref{eq:M}'}
\end{equation}
Multiplying, for all $\dg=0:2d$, the equation of the system 
(\ref{eq:M}')
corresponding to $\dg$ with $\binom{d+1+\dg}{\dg}$
and noting that, by convention, $\binom{d+1+\dg}{\dg} =1$ for $\dg=0$,
we can write the system
(\ref{eq:M}')
in the simpler equivalent form
$M_0 \mc = e_1$.
Using the fact from \cite[Sect.~8]{deboor} 
that the $k$th divided difference is linked to the
the $k$th derivative by 
$\divdif{\tsft} t^{d+1+\dg} 
= (d+1+\dg)\cdots(1+\dg)\ \xi_\sft^\dg$,
for distinct $\xi_\sft \in (t_\sft..t_{\sft+d+1}]$, 
we see that $M_0$ is a 
Vandermonde matrix, hence invertible\footnote{
Mirzaee et al.\ \cite[p.90]{journals/jscic/MirzaeeRK12}
state that systems of type (\ref{eq:convol}') are non-singular
and refer to the exposition in \cite{journals/moc/CockburnLSS03} for a
a proof of existence and uniqueness of the solution.
I was unable to spot it there.}. That is $\mc := M^{-1}_0\be_1$ 
is well-defined.

To show that $\mc$ solves the polynomial system
(\ref{eq:convol}') for all $x$, not just $x=0$, 
we define the $2d+1$ linearly independent functionals
$F_k$, $k=0,\ldots,2d$.
The functional $F_k$ differentiates (each entry of $M(x)$)
$k$-times with respect to $x$ and then evaluates at zero.
Applying $F_k$ to both sides of (\ref{eq:convol}'') yields the 
system  
\begin{align}
   M_k \mc
   = \binom{d-k+\dg}{d+1} e_{k+1}, \qquad 
   M_k(j,:) :=  
   \begin{cases}
      M_0(j-k,:), & \text{ if } j>k,\\
      0, & \text{ else.}
   \end{cases}
   \tag{\ref{eq:convol}$^k$}
\end{align}
(Here $M_k(j,:)$ denotes the $j$th row of the matrix $M_k$.)
Since $M_k\mc =  M_k M^{-1}_0 \be_1 = \be_{k+1}$ and 
$\binom{d-k+\dg}{d+1}=1$ for $\dg=k+1$, we see that
the choice $\mc := M^{-1}_0 \be_1$ satisfies
all $2d+1$ systems of equations (\ref{eq:convol}$^k$). 
This implies that the system of polynomial equations (\ref{eq:convol}$'$) 
is satisfied by $\mc$.
%
\end{proof}


If the knots are strictly increasing,
we can expand the divided differences of polynomials to make
the expression for $M_0$ more explicit.

\begin{corollary} [single knots]
If, for $i=-d:2d$, $t_i<t_{i+1}$ then
\begin{align}
   M_0 = 
   \left[ \begin{matrix}
   \sum^{\sft+d+1}_{\ell=\sft} 
   \frac{(x-x_\ell)^{d+1+\dg}}
   {\prod^{\sft+d+1}_{j=\sft,j\ne \ell} (x_j-x_\ell)}
   \end{matrix} \right]_{\dg=0:2d, \sft = -d:d}.
\end{align}
\end{corollary}

When the knot spacing is uniform,
i.e.\ $t_{i+1}-t_i = t_i-t_{i-1}$, then it is good to 
\emph{symmetrize} the construction about $t=0$. 
That is, we define the knot sequence to be
$\tau := [-d-\gs : d-\gs]$, $\gs := \frac{d+1}{2}$.
For example, for $d=1$, $\tau = [-2,-1,0]$
and for $d=2$, $\tau = [-3.5,-2.5,-1.5,-.5,0.5]$.
Then for $d=2$, the first knot subsequence is $[-3.5:-0.5]$ and,
symmatrically, the last knot subsequence is $[0.5:3.5]$.

\begin{corollary}[uniform knots]
For uniform knots 
\begin{equation}
   M_0
   =
   \left[ \begin{matrix}
   \frac{1}{(d+1)!}
   \sum^{d+1}_{\ell=0}
   (-1)^\ell \binom{d+1}{\ell} (\sft+\ell)^{d+1+\dg}
   \end{matrix} \right]_{\dg=0:2d, \sft \in \tau}.
   \label{eq:tkk}
\end{equation}
\end{corollary}
For comparison,
Mirzagar \cite{mirz:2014} characterizes the
symmetric SIAC kernel coefficients by the relations
$\tilde{M}(x) \mc = [x^\dg]_{\dg=0:2d}$, where
\begin{align}
   \tilde{M}(x)
   := \left[ \begin{matrix}
   \sum^\dg_{\ell=0}  (-1)^\ell \sft^\ell
   \binom{\dg}{\ell} 
   (\bsp(\cdot|0:2d)\convol (\cdot)^{\dg-\ell})(x)
   \end{matrix} \right]_{\dg=0:2d, \sft =-d:d}.
\end{align}



\noindent
The explicit form \eqref{eq:tkk} makes it easy to confirm
a conjecture by Kirby and Ryan that
the optimal SIAC coeffients in the uniform case
are rational numbers:
by Cramer's rule, 
\begin{align}
   \mc &= \frac{\det[\be_1 M_0(:,2:2d+1)]}{\det M_0}
\end{align}
and the determinants only contain rational numbers.
For example, the optimal symmetric SIAC spline convolution coefficients
for degree $d$ are (omitted entries in slots $d+2:2d+1$
indicated by ``$\ldots$'' are defined by symmetry):
\begin{align*}
   d &=1 :\quad [-1, 14, -1]/12,
   \\
   d &=2 :\quad [-37, 388, -2622, 388, -37]/1920,
   \\
   d &=3 :\quad [-82, 933, -5514, 24446, -5514, 933, -82]/
      15120,
      \\
   d &=4 :\quad [-153617, 1983016, -12615836, 54427672, -180179750,\ldots]/
       92897280,
       \\
   d &=5 :\quad [-4201, 61546, -437073, 2034000, -7077894, 18830604,\ldots]/
       7983360.
\end{align*}

\section{Conclusion}
Especially in the presence of non-uniform knots, where pre-tabulation may 
not be practical, it is good to have explicit formulas for the entries of
the SIAC coefficient matrix $M_0$ in terms of divided differences.
Numerically stable implementations of divided differences are
well-known.
The paper may serve as a building block towards addressing the important 
issue of choosing good non-uniform knot sequences for SIAC post-processing.
The author conjectures that the subtle issue of knot selection 
is also closely related to a classical theorem of spline theory.

\paragraph{Acknowledgement}
This work was supported in part by NSF grant CCF-1117695 and
NIH R01 LM011300-01. 
Jennifer Ryan, Mahsa Mirzagar, Mike Kirby and Xiaozhou Li introduced
me to the topic and patiently answered my questions.
\bibliographystyle{alpha}
\bibliography{p}

\begin{thebibliography}{WRKH09}

\bibitem[BS77]{Bramble:1977:HOL}
J.~H. Bramble and A.~H. Schatz.
\newblock Higher order local accuracy by averaging in the finite element
  method.
\newblock {\em Mathematics of Computation}, 31(137):94--111, January 1977.

\bibitem[CKRS07]{journals/siamsc/CurtisKRS07}
Sean Curtis, Robert~M. Kirby, Jennifer~K. Ryan, and Chi-Wang Shu.
\newblock Postprocessing for the discontinuous {G}alerkin method over
  nonuniform meshes.
\newblock {\em SIAM J. Scientific Computing}, 30(1):272--289, 2007.

\bibitem[CLSS03]{journals/moc/CockburnLSS03}
Bernardo Cockburn, Mitchell Luskin, Chi-Wang Shu, and Endre S{\"u}li.
\newblock Enhanced accuracy by post-processing for finite element methods for
  hyperbolic equations.
\newblock {\em Math. Comput}, 72(242):577--606, 2003.

\bibitem[dB02]{Boor:2002:BB}
C.W. de~Boor.
\newblock B-spline {B}asics.
\newblock In M.~Kim G.~Farin, J.~Hoschek, editor, {\em Handbook of Computer
  Aided Geometric Design}. Elsevier, 2002.

\bibitem[dB05]{deboor}
C.~W. de~Boor.
\newblock Divided differences.
\newblock {\em Surveys in Approximation Theory}, 1:46--69, 2005.

\bibitem[MJRK11]{journals/siamnum/MirzaeeJRK11}
Hanieh Mirzaee, Liangyue Ji, Jennifer~K. Ryan, and Robert~M. Kirby.
\newblock Smoothness-increasing accuracy-conserving ({SIAC}) postprocessing for
  discontinuous {G}alerkin solutions over structured triangular meshes.
\newblock {\em SIAM J. Numerical Analysis}, 49(5):1899--1920, 2011.

\bibitem[MKRK13]{journals/siamsc/MirzaeeKRK13}
Hanieh Mirzaee, James King, Jennifer~K. Ryan, and Robert~M. Kirby.
\newblock Smoothness-increasing accuracy-conserving filters for discontinuous
  {G}alerkin solutions over unstructured triangular meshes.
\newblock {\em SIAM J. Scientific Computing}, 35(1), 2013.

\bibitem[ML78]{MockLax}
M.S. Mock and P.D. Lax.
\newblock The computation of discontinuous solutions of linear hyperbolic
  equations.
\newblock {\em Comm. Pure Appl. Math.}, 31(4):423--430, 1978.

\bibitem[MR14]{mirz:2014}
Mahsa Mirzargar and Jennifer~K. Ryan.
\newblock Filtering: A unified view from approximation theory to the
  post-processing of discontinuous {G}alerkin solutions, 2014.
\newblock Mini-symposium at Icosahom 2014, Salt Lake City, July 2014.

\bibitem[MRK12]{journals/jscic/MirzaeeRK12}
Hanieh Mirzaee, Jennifer~K. Ryan, and Robert~M. Kirby.
\newblock Efficient implementation of smoothness-increasing accuracy-conserving
  ({SIAC}) filters for discontinuous {G}alerkin solutions.
\newblock {\em J. Sci. Comput}, 52(1):85--112, 2012.

\bibitem[MRK14]{journals/jscic/MirzaeeRK14}
Hanieh Mirzaee, Jennifer~K. Ryan, and Robert~M. Kirby.
\newblock Smoothness-increasing accuracy-conserving ({SIAC}) filters for
  discontinuous {G}alerkin solutions: Application to structured tetrahedral
  meshes.
\newblock {\em J. Sci. Comput}, 58(3):690--704, 2014.

\bibitem[Sch81]{Schumaker:1981:SF}
L.L. Schumaker.
\newblock {\em Spline Functions: Basic Theory}.
\newblock 1981.
\newblock Wiley, New York.

\bibitem[WRKH09]{journals/jscic/WalfischRKH09}
David Walfisch, Jennifer~K. Ryan, Robert~M. Kirby, and Robert Haimes.
\newblock One-sided smoothness-increasing accuracy-conserving filtering for
  enhanced streamline integration through discontinuous fields.
\newblock {\em J. Sci. Comput}, 38(2):164--184, 2009.

\end{thebibliography}
\end{document}